# *Fuzzy derivations KU-ideals on KU-algebras*

## BY

*Samy M.Mostafa [1]*  ,  *Ahmed Abd-eldayem[2]*

samymostafa@yahoo.com      ahmedabdeldayem88@yahoo.com

1,2Department of mathematics -Faculty of Education -Ain Shams University Roxy, Cairo, Egypt

*Abstract*. In this manuscript, we introduce a new concept, which is called fuzzy left (right) derivations KU- ideals in KU-algebras. We state and prove some theorems about fundamental properties of it. Moreover, we give the concepts of the image and the pre-image of fuzzy left (right) derivations KU-ideals under homomorphism of KU- algebras and investigated some its properties. Further, we have proved that every the image and the pre-image of fuzzy left (right) derivations KU-ideals under homomorphism of KU- algebras are fuzzy left (right) derivations KU-ideals. Furthermore, we give the concept of the Cartesian product of fuzzy left (right) derivations KU - ideals in Cartesian product of KU – algebras.

## *1. Introduction*

As it is well known, BCK and BCI-algebras are two classes of algebras of logic. They were introduced by Imai and Iseki [10,11,12] and have been extensively investigated by many researchers. It is known that the class of BCK-algebras is a proper sub class of the BCI-algebras.The class of all *BCK*-algebras is a quasivariety. Is´*e*ki posed an interesting problem (solved by Wro´*n*ski [24]) whether the class of *BCK*-algebras is a variety. In connection with this problem, Komori [15] introduced a notion of *BCC*-algebras, and Dudek [7] redefined the notion of *BCC*-algebras by using a dual form of the ordinary definition in the sense of Komori. Dudek and Zhang [8] introduced a new notion of ideals in *BCC*-algebras and described connections between such ideals and congruences .

C.Prabpayak and U.Leerawat ( [22 ], [23 ]) introduced a new algebraic structure which is called KU - algebra . They gave the concept of homomorphisms of KU- algebras and investigated some related properties. Several authors [2,3 ,5 ,6,9,14] have studied derivations in rings and near rings. Jun and Xin [13] applied the notion of derivations in ring and near-ring theory to *BCI*-algebras, and they also introduced a new concept called a regular derivation in *BCI* -algebras. They investigated some of its properties, defined a *d* -derivation ideal and gave conditions for an ideal to be *d*-derivation. Later, Hamza and Al-Shehri [1], defined a left derivation in *BCI*-algebras and investigated a regular left derivation. Zhan and Liu [27 ] studied f-derivations in BCI-algebras and proved some results. G. Muhiuddin etl [20,21] introduced the notion of $(\alpha,\beta)$-derivation in a BCI-algebra and investigated related properties. They provided a condition for a $(\alpha,\beta)$- derivation to be regular. They also introduced the concepts of a $d_{(\alpha,\beta)}$-invariant $(\alpha,\beta)$ -derivation and α-ideal, and then they investigated their relations. Furthermore, they obtained some results on regular $(\alpha,\beta)$ - derivations. Moreover, they studied the notion of *t*-derivations on BCI-algebras and obtained some of its related properties. Further, they characterized the notion of *p*-semi-simple BCI-algebra *X* by using the notion of *t*-derivation. Later, Mostafa et al [18,19], introduced the notions of $(\ell,r)$ -$((r,\ell))$ -derivation of a KU-algebra and some related properties are explored. The concept of fuzzy sets was introduced by Zadeh [26]. In 1991, Xi [25] applied the concept of fuzzy sets to BCI, BCK, MV-algebras .Since its inception, the theory of fuzzy sets , ideal theory and its fuzzification has developed in many directions and is finding applications in a wide variety of fields. Mostafa et al, in 2011[17] introduced the notion of fuzzy KU-ideals of KU-algebras and then they investigated several basic properties which are related to fuzzy KU-ideals.

In this paper, we introduce the notion of fuzzy left (right) derivations KU- ideals in KU - algebras. The homomorphic image ( preimage) of fuzzy left (right) - derivations KU- ideals in KU - algebras under homomorprhism of a *KU*-algebras are discussed . Many related results have been derived.

## 2. Preliminaries

In this section, we recall some basic definitions and results that are needed for our work.

*Definition 2.1 [22,23]* Let $X$ be a set with a binary operation $*$ and a constant 0. $(X, *, 0)$ is called KU-algebra if the following axioms hold : $\forall\, x, y, z \in X$ :

$(KU_1)$  $(x * y) * [(y * z) * (x * z)] = 0$

$(KU_2)$  $x * 0 = 0$

$(KU_3)$  $0 * x = x$

$(KU_4)$  if $x * y = 0 = y * x$ implies $x = y$

Define a binary relation $\leq$ by : $x \leq y \Leftrightarrow y * x = 0$, we can prove that $(X, \leq)$ *is a partially ordered set*.

Throughout this article, $X$ will denote a KU-algebra unless otherwise mentioned

*Corollary 2.2 [17,22]* In KU-algebra the following identities are true for all $x, y, z \in X$ :

(i) $z * z = 0$

(ii) $z * (x * z) = 0$

(iii) If $x \leq y$ implies that $y * z \leq x * z$

(v) $z * (y * x) = y * (z * x)$

(vi) $y * [(y * x) * x] = 0$

*Definition 2.3 [22,23]* A subset S of KU-algebra X is called sub algebra of X if $x * y \in S$, whenever $x, y \in S$

*Definition 2.4 [22,23]* A non empty subset $A$ of KU-algebra X is called ideal of X if it is satisfied the following conditions:

(i) $0 \in A$

(ii) $y * z \in A$, $y \in A$   implies   $z \in A$   $\forall y, z \in X$.

***Definition 2.5 [17]*** A non-empty subset A of a KU-algebra X is called an KU ideal of X if it satisfies the following conditions:

(1) $0 \in A$,

(2) $x * (y * z) \in A$, $y \in A$ implies $x * z \in A$, for all $x, y, z \in X$

***Definition 2.6 [17]*** Let X be a KU-algebra, a fuzzy set μ in X is called fuzzy sub-algebra if it satisfies:

(S$_1$)  $\mu(0) \geq \mu(x)$,

(S$_2$)  $\mu(x) \geq \{\mu(x * y), \mu(y)\}$ for all $x, y \in X$.

***Definition 2.7 [17]*** Let X be a KU-algebra, a fuzzy set μ in X is called a fuzzy KU-ideal of X if it satisfies the following conditions:

(F$_1$)  $\mu(0) \geq \mu(x)$,

(F$_2$)  $\mu(x * z) \geq \min\{\mu(x * (y * z)), \mu(y)\}$.

***Definition 2.8*** For elements $x$ and $y$ of KU-algebra ($X, *, 0$), we denote
$x \wedge y = (x * y) * y$.

***Definition 2.9 [18]*** Let $X$ be a KU-algebra. A self map $d : X \to X$ is a left–right derivation (briefly, $(\ell, r)$-derivation) of $X$ if it satisfies the identity

$$d(x * y) = (d(x) * y) \wedge (x * d(y)) \ \forall x, y \in X$$

If d satisfies the identity

$$d(x * y) = (x * d(y)) \wedge (d(x) * y) \ \forall x, y \in X$$

d is called right-left derivation (briefly, $(r, \ell)$-derivation) of $X$. Moreover, if $d$ is both $(\ell, r)$ and $(r, \ell) - derivation$ then $d$ is called a derivation of $X$.

***Definition 2.10 [18]*** A derivation of KU-algebra is said to be regular if $d(0) = 0$.

*Lemma 2.11[18]* A derivation $d$ of KU-algebra $X$ is regular.

*Example 2.12 [18]* Let $X = \{0,1,2,3,4\}$ be a set in which the operation $*$ is defined as follows:

| * | 0 | 1 | 2 | 3 | 4 |
|---|---|---|---|---|---|
| 0 | 0 | 1 | 2 | 3 | 4 |
| 1 | 0 | 0 | 2 | 2 | 4 |
| 2 | 0 | 0 | 0 | 1 | 4 |
| 3 | 0 | 0 | 0 | 0 | 4 |
| 4 | 0 | 1 | 1 | 1 | 0 |

Using the algorithms in Appendix A, we can prove that (X, *, 0) is a KU-algebra. Define a map $d: X \to X$ by

$$d(x) = \begin{cases} 0 & \text{if } x = 0,1,2,3 \\ 4 & \text{if } x = 4 \end{cases}$$

Then it is easy to show that $d$ is both a $(\ell, r)$ and $(r, \ell)$-derivation of $X$.

Example 2.12. Let $Z^+ \cup \{0\}$ be the set of all positive integers and 0. The operation $(*)$ on $Z^+ \cup \{0\}$ is defined as follows: $x*y = y-x$, where "−" the minus operation. Define a binary relation $\leq$ on $Z^+ \cup \{0\}$ by: $x \leq y \Leftrightarrow y*x = 0$. Then $X = (Z^+ \cup \{0\}, *, 0)$ is a KU-algebra. We define a map $d: X \to X$ by $d(x) = x - 1$ for all $x \in Z^+$. Then $\forall x, y \in X$, we have

d(x*y)=d(y−x)=y−x−1………………………………………………………….(I),

d(x)*y = y − d(x)=y− (x−1) =1+ y − x and x*d(y) =d(y) −x =y−1−x=y−x−1, but

d(x)*y∧ x*d(y) =((1+y −x )* (y−x−1))* (y−x−1)= (y−x−1) −[(y−x−1) −(1+y−x )]=

= y+1− x……………………………. (II)

From (I) and (II), d is not $(\ell, r)$ derivation of X.
On other hand

x*d(y)= d(y) −x=y−x−1 ,   d(x)*y= y− d(x) =y−(x−1)=y+1−x , but

x*d(y) ∧ d(x)*y = [ (x*d(y))* (d(x)*y)]* (d(x)*y)= ( d(x) −y) −[(y −d(x)) − (d(y) −x)]

$$= y-x-1 \dots\dots\dots\dots\dots\dots\dots\dots\dots \text{(III)}$$

From (I) and (III), d is $(r, \ell)$ derivation of X. Hence $(r, \ell)$-derivation and $(\ell, r)$ derivation are not coincide.

***Proposition 2.13[18]*** Let $X$ be a KU-algebra with partial order $\leq$, and let $d$ be a derivation of $X$. Then the following hold $\forall\, x, y \in X$ :

**(i)** $d(x) \leq x$.

**(ii)** $d(x * y) \leq d(x) * y$.

**(iii)** $d(x * y) \leq x * d(y)$.

**(v)** $d(x * d(x)) = 0$.

**(vi)** $d^{-1}(0) = \{x \in X \mid d(x) = 0\}$ is a sub algebra of $X$.

***Definition 2.14 [18]*** Let $X$ be a KU-algebra and $d$ be a derivation of $X$.
Denote $Fix_d(X) = \{x \in X : d(x) = x\}$.

***Proposition 2.15[18]*** Let $X$ be a KU-algebra and $d$ be a derivation of $X$. Then $Fix_d(X)$ is a sub algebra of $X$.

## 3. *Fuzzy derivations KU- ideals of KU-algebras*

In this section, we will discuss and investigate a new notion called fuzzy- left derivations KU - ideals of KU - algebras and study several basic properties which are related to fuzzy left derivations KU - ideals.

***Definition 3.1*** Let X be a KU-algebra and $d : X \to X$ be self map. A non - empty subset A of a KU-algebra X is called left derivations KU ideal of X if it satisfies the following conditions:

(1) $0 \in A$,

(2) $d(x) * (y * z) \in A$, $d(y) \in A$ implies $d(x * z) \in A$, for all $x, y, z \in X$

***Definition 3.2*** Let X be a KU-algebra and $d : X \to X$ be self map. A non - empty subset A of a KU-algebra X is called right derivations KU ideal of X if it satisfies the following conditions:

(1) $0 \in A$,

(2) $x * d(y * z) \in A$, $d(y) \in A$ implies $d(x * z) \in A$, for all $x, y, z \in X$.

***Definition 3.3*** Let X be a KU-algebra and $d : X \to X$ be self map. A non - empty subset A of a KU-algebra X is called derivations KU -ideal of X if it satisfies the following conditions:

(1) $0 \in A$,

(2) $d(x * (y * z)) \in A$, $d(y) \in A$ implies $d(x * z) \in A$, for all $x, y, z \in X$

***Definition 3.4*** Let X be a KU-algebra and $d : X \to X$ be self map. A fuzzy set $\mu : X \to [0,1]$ in X is called a fuzzy left derivations KU-ideal(briefly, $(F, \ell)$-derivation) of X if it satisfies the following conditions:

(F$_1$)   $\mu(0) \geq \mu(x)$,

(FL$_2$)   $\mu(d(x * z)) \geq \min\{\mu(d(x)*(y*z)), \mu(d(y))\}$.

***Definition 3.5*** Let X be a KU-algebra and $d: X \to X$ be self map. A fuzzy set $\mu: X \to [0,1]$ in X is called a fuzzy right derivations KU-ideal(briefly, $(F,r)$-derivation) of X if it satisfies the following conditions:

($F_1$)  $\mu(0) \geq \mu(x)$.

($FR_2$)  $\mu(d(x*z)) \geq \min\{\mu(x*d(y*z)), \mu(d(y))\}$.

***Definition 3.6*** Let X be a KU-algebra and $d: X \to X$ be self map. A fuzzy set $\mu: X \to [0,1]$ in X is called a fuzzy derivations KU-ideal, if it satisfies the following conditions

($F_1$)  $\mu(0) \geq \mu(x)$.

($F_2$)  $\mu(d(x*z)) \geq \min\{\mu(d(x*(y*z))), \mu(d(y))\}$.

***Remark 3.7 (I)*** If d is fixed, the definitions (3.1, 3.2, 3.3) gives the definition KU-ideal.

 **(II)** If d is fixed, the definitions (3.4, 3.5, 3.6) gives the definition fuzzy KU-ideal.

***Example 3.8*** Let $X = \{0,1,2,3,4\}$ be a set in which the operation * is defined as follows:
Using the algorithms in Appendix A, we can prove that (X, *, 0) is a KU-algebra.

| * | 0 | 1 | 2 | 3 | 4 |
|---|---|---|---|---|---|
| 0 | 0 | 1 | 2 | 3 | 4 |
| 1 | 0 | 0 | 2 | 3 | 3 |
| 2 | 0 | 1 | 0 | 1 | 4 |
| 3 | 0 | 0 | 0 | 0 | 3 |
| 4 | 0 | 0 | 0 | 0 | 0 |

Define a self map $d: X \to X$ by

$$d(x) = \begin{cases} 0 & \text{if } x = 0,1,2,3 \\ 4 & \text{if } x = 4 \end{cases}.$$

Define a fuzzy set $\mu : d(X) \to [0,1]$, by $d(\mu(0)) = t_0$, $\mu(d(1)) = \mu(d(2)) = t_1$, $\mu(d(3)) = \mu(d(4)) = t_2$, where $t_0, t_1, t_2 \in [0,1]$ with $t_0 > t_1 > t_2$. Routine calculations give that $\mu$ is a fuzzy left (right)- derivations KU- ideal of KU- algebra X.

**Lemma 3.9** Let $\mu$ be a fuzzy left derivations KU - ideal of KU - algebra X, if the inequality $x*y \leq d(z)$ holds in X, then $\mu(d(y)) \geq \min\{\mu(d(x)), \mu(z)\}$.

Proof. Assume that the inequality $x*y \leq d(z)$ holds in X, then

$d(z) * (x * y) = 0$, $(z) * (x * y) = 0$, since $d(z) \leq z$ from (Proposition 2.13(i)) and by (FL$_2$), we have $\mu(d(z*y)) \geq \min\{\mu(d(z)*(x*y)), \mu(d(x))\} = \min\{\mu(0), \mu(d(x))\} = \mu(d(x))$

Put $z=0$, we have $\mu(d(0*y)) = \mu(d(y)) \geq \min\{\mu(x*y), \mu(d(x))\}$ ………… (i),

but $\mu(x*y) \geq \min\{\mu(x*(z*y)), \mu(z)\} = \min\{\mu(z*(x*y)), \mu(z)\}$

$$= \min\{\mu(0), \mu(z)\} = \mu(z) \quad \ldots\ldots\ldots(ii)$$

From (i), (ii), we get $\mu(d(y)) \geq \min\{\mu(z), \mu(d(x))\}$, this completes the proof.

**Lemma 3.10** If $\mu$ is a fuzzy left derivations KU - ideal of KU - algebra X and if $x \leq d(y)$, then $\mu(d(x)) \geq \mu(d(y))$.

Proof. If $x \leq d(y)$, then $d(y) * x = 0$, $y*x=0$ since $d(y) \leq y$ (from Proposition 2.13(i)) this together with $0 * x = x$ and $\mu(0) \geq \mu(y)$, we get

$\mu(d(0*x)) = \mu(d(x)) \geq \min\{\mu(d(0)*(y*x)), \mu(d(y))\} = \min\{\mu(0*0), \mu(d(y))\} =$
$$= \min\{\mu(0), \mu(d(y))\} = \mu(d(y)).$$

**Proposition 3.11** The intersection of any set of fuzzy left derivations KU - ideals of KU – algebra X is also fuzzy left derivations KU - ideal.

Proof. let $\{\mu_i\}$ be a family of fuzzy left derivations KU - ideals of KU- algebra X, then for any $x, y, z \in X$,

$(\bigcap \mu_i)(0) = \inf(\mu_i(0)) \geq \inf(\mu_i(d(x))) = (\bigcap \mu_i)(d(x))$ and

$(\bigcap \mu_i)(d(x*z)) = \inf(\mu_i(d(x*z)) \geq \inf(\min\{\mu_i((d(x)*(y*z)), \mu_i(d(y))\} =$
$\min\{\inf(\mu_i((d(x)*(y*z))), \inf(\mu_i(d(y))\} = \min\{(\bigcap \mu_i)((d(x)*(y*z))), (\bigcap \mu_i)(d(y))\}.$

This completes the proof.

***Lemma 3.12*** The intersection of any set of fuzzy right derivations KU - ideals of KU – algebra X is also fuzzy right derivations KU - ideal.

***proof.*** Clear

***Theorem 3.13*** Let µ be a fuzzy set in X then µ is a fuzzy left derivations KU- ideal of X if and only if it satisfies :

For all α∈[0,1]),U (µ , α) ≠ φ implies U(µ ,α) is KU- ideal of X……… (A),

where U (µ , α) = {x ∈ X / µ (d(x)) ≥ α} .

Proof . Assume that µ is a fuzzy left derivations KU- ideal of X , let α ∈ [0 , 1] be such that U (µ , α) ≠ φ , and x , y ∈ X such that x ∈ U (µ , α) , then µ (d(x)) ≥ α and so by (FL$_2$) , µ (d( y * 0)) = µ (d(0)) ≥ min { µ ( d(y )* (x * 0) ) , µ (d(x))}=
min{µ (d(y) * 0), µ (d(x))} = min {µ (0) , µ (d(x))} = α , hence  0 ∈ U (µ , α) .
Let  d(x) * (y * z) ∈ U (µ, α ) , d(y) ∈ U (µ, α), It follows from(FL$_2$) that
µ (d(x * z)) ≥ min {µ (d(x) * (y * z)) , µ (d(y))} = α , so that x * z ∈ U (µ, α) .
Hence U (µ, α ) is KU - ideal of X .
Conversely, suppose that µ satisfies (A) , let x , y , z ∈ X be such that
µ (d(x * z)) < min {µ (d(x) * (y * z)) , µ (d(y))},taking
β0 = 1/2 {µ (d(x * z)) + min {µ (d(x) * (y * z)) , µ(d(y)) } , we have
β0 ∈ [0,1] and µ ( d(x * z)) < β0 < min {µ (d(x) * (y * z)) , µ(d(y)) } ,it follows that
d(x) * (y * z) ∈ U (µ, β0) and d(x * z) ∈ U (µ, β0) , this is a contradiction and therefore µ is a fuzzy left derivations KU - ideal of X .

***Theorem 3.14*** Let µ be a fuzzy set in X then µ is a fuzzy right derivations KU- ideal of X if and only if it satisfies : For all α∈[0,1]),U (µ , α) ≠ φ implies U(µ ,α) is KU- ideal of X.

***Proposition 3.15*** If µ is a fuzzy left derivations KU - ideal of X , then

$$\mu (d(x) * (x * y)) \geq \mu (d(y))$$

proof. Taking z = x * y in (FL2) and using (ku2) and (F1) , we get

μ(d(x) * (x * y)) ≥ min { μ (d(x) * (y * (x * y)) , μ(d(y)) } = min {μ(d(x) * (x * (y * y)) , μ(d(y)) } = min {μ( d(x) * (x * 0)) , μ(d(y)) }= min {μ (0) , μ (d(y)) } = μ (d(y)).

*Definition3.16* Let μ be a fuzzy left derivations KU - ideal of KU - algebra X ,.the KU - ideals $\mu_t$ , t∈ [0,1] are called level KU - ideal of μ .

*Corollary3.17* Let I be an KU - ideal of KU - algebra X , then for any fixed number t in an open interval (0,1) , there exist a fuzzy left derivations  KU – ideal μ of X such that $\mu_t$ = I .

proof. The proof is similar the corollary 4.4 [17] .

# 4 *Image (Pre-image) of fuzzy derivations  KU-ideals under homomorphism*

In this section, we introduce the concepts of the image and the pre-image of fuzzy left derivations  KU-ideals in KU-algebras under homomorphism.

*Definition 4.1* Let f be a mapping from the set X to a set Y. If μ is a fuzzy subset of X, then the fuzzy subset β of Y defined by

$$f(\mu)(y) = \beta(y) = \begin{cases} \sup_{x \in f^{-1}(y)} \mu(x), & \text{if } f^{-1}(y) = \{x \in X, f(x) = y\} \neq \phi \\ 0 & \text{otherwise} \end{cases}$$

is said to be the image of μ under f.

Similarly if β is a fuzzy subset of Y , then the fuzzy subset μ = β o f in X ( i.e the fuzzy subset defined by μ (x) = β (f (x)) for all x ∈ X ) is called the primage of β under f .

***Theorem 4.2*** An onto homomorphic preimage of a fuzzy left derivations KU - ideal is also a fuzzy left derivations KU - ideal .

Proof. Let $f : X \to X`$ be an onto homomorphism of KU - algebras , $\beta$ a fuzzy left derivations KU - ideal of $X`$ and $\mu$ the preimage of $\beta$ under f , then $\beta$ (f (d(x)) = $\mu$ (d(x)) , for all $x \in X$ . Let $x \in X$ , then $\mu$ (d(0)) = $\beta$ (f(d(0))) $\geq \beta$ (f (d(x))) = $\mu$ (d(x)).

Now let x , y , z $\in X$ , then $\mu$ (d(x $*$ z)) = $\beta$ (f (d(x $*$ z))) $\geq$

$$\min \{\beta(f (d(x)) * (f (y) * f(z)), \beta(f (d(y)))\} =$$

$$\min \{ \beta (f(d(x))*(y * z)), \beta (f (d(y)) \} =$$

$$\min \{\mu(d(x) * (y * z))) , \mu(d (y))\} .$$

The proof is completed.

***Definition 4.3* [4]** A fuzzy subset $\mu$ of X has sup property if for any subset T of X , there exist $t_0 \in T$ such that , $\mu(t_0) = \underset{t \in T}{SUP} \mu(t)$.

***Theorem 4.4*** Let $f : X \to Y$ be a homomorphism between KU - algebras X and Y .

For every fuzzy left derivations KU - ideal $\mu$ in X , f ($\mu$) is a fuzzy left derivations KU - ideal of Y .

Proof. By definition $\beta(d(y')) = f(\mu)(d(y')) = \underset{d(x) \in f^{-1}((d(y'))}{\sup} \mu(d(x))$ for all $y' \in Y$ and $\sup \phi = 0$. We have to prove that $\beta(d(x' * z')) \geq \min\{\beta(d(x')*(y'*z'), \beta(d(y'))\}$,

$\forall$ x`, y`, z`$\in Y$. Let $f : X \to Y$ be an onto a homomorphism of KU - algebras , $\mu$ a fuzzy left derivations KU - ideal of X with sup property and $\beta$ the image of $\mu$ under f , since $\mu$ is a fuzzy left derivations KU - ideal of X , we have $\mu(d(0)) \geq \mu(d(x))$ for all $x \in X$ . Note that $0 \in f^{-1}(0`)$ , where 0 , 0` are the zero of X and Y respectively

Thus, $\beta(d(0')) = \sup_{d(t) \in f^{-1}(d(0'))} \mu(d(t)) = \mu(d(0)) = \mu(0) \geq \mu(d(x))$, for all $x \in X$,

which implies that $\beta(d(0')) \geq \sup_{d(t) \in f^{-1}(d(x'))} \mu(d(t)) = \beta(d(x'))$, for any $x' \in Y$. For any $x', y', z' \in Y$, let $d(x_0) \in f^{-1}(d(x'))$, $d(y_0) \in f^{-1}(d(y'))$, $d(z_0) \in f^{-1}(d(z'))$ be Such that $\mu(d(x_0 * z_0)) = \sup_{d(t) \in f^{-1}(d(x'*z'))} \mu(d(t))$, $\mu(y_0) = \sup_{d(t) \in f^{-1}(d(y'))} \mu(d(t))$

and

$$\mu(d(x_0) * (y_0 * z_0)) = \beta\{f(d(x_0) * (y_0 * z_0))\} = \beta(d(x') * (y' * z')) =$$
$$\sup_{(d((x_0)*(y_0*z_0)) \in f^{-1}(d(x')*(y'*z'))} \mu(d(x_0) * (y_0 * z_0))\} = \sup_{d(t) \in f^{-1}(d(x')*(y'*z'))} \mu(d(t)) \cdot$$

Then

$$\beta(d(x' * z')) = \sup_{d(t) \in f^{-1}(d(x'*z'))} \mu(d(t)) = \mu(d(x_0 * z_0)) \geq \min\{\mu(d(x_0) * (y_0 * z_0)), \mu(d(y_0))\} =$$

$$\min\left\{\sup_{d(t) \in f^{-1}(d(x')*(y'*z'))} \mu(d(t)), \sup_{d(t) \in f^{-1}(d(y'))} \mu(d(t))\right\} = \min\{\beta(d(x') * (y' * z')), \beta(d(y'))\}.$$

Hence $\beta$ is a fuzzy left derivations KU-ideal of Y.

***Theorem 4.5*** An onto homomorphic preimage of a fuzzy right derivations KU - ideal is also a fuzzy right derivations KU - ideal

***Theorem 4.6*** Let $f : X \to Y$ be a homomorphism between KU - algebras X and Y.

For every fuzzy right derivations KU - ideal μ in X , f (μ) is a fuzzy right derivations KU - ideal of Y.

proof. Clear

## 5. Cartesian product of fuzzy left derivations KU-ideals

*Definition 5.1*[4] A fuzzy µ is called a fuzzy relation on any set S, if µ is a fuzzy subset µ : $S \times S \to [0,1]$.

*Definition 5.2* [4] If µ is a fuzzy relation on a set S and β is a fuzzy subset of S, then µ is a fuzzy relation on β if µ (x, y) ≤ min {β (x), β (y)}, ∀ x, y ∈ S.

*Definition 5.3 [4]* Let µ and β be fuzzy subset of a set S, the Cartesian product of µ and β is define by (µ × β) (x, y) = min {µ (x), β (y)}, ∀ x, y ∈ S.

*Lemma 5.4[4]* Let µ and β be fuzzy subset of a set S, then

(i) $\mu \times \beta$ is a fuzzy relation on S.

(ii) $(\mu \times \beta)_t = \mu_t \times \beta_t$ for all t ∈ [0,1].

*Definition 5.5* If µ is a fuzzy left derivations relation on a set S and β is a fuzzy left derivations subset of S, then µ is a fuzzy left derivations relation on β if

µ(d (x, y)) ≤ min {β(d (x)), β(d (y))}, ∀ x, y ∈ S.

*Definition 5.6 [4]* Let µ and β be fuzzy left derivations subset of a set S, the Cartesian product of µ and β is define by (µ × β)(d (x, y)) = min {µ(d (x)), β(d (y))}, ∀ x, y ∈ S

*Lemma 5.7[4]* Let µ and β be fuzzy subset of a set S, then

(i) $\mu \times \beta$ is a fuzzy relation on S,

(ii) $(\mu \times \beta)_t = \mu_t \times \beta_t$ for all t ∈ [0,1].

***Definition 5.8*** If β is a fuzzy left derivations subset of a set S , the strongest fuzzy relation on S , that is a fuzzy derivations relation on β is $\mu_\beta$ given by

$\mu_\beta(d(x, y)) = \min \{\beta(d(x)), \beta(d(y))\}, \forall\ x, y \in S$ .

***Lemma 5.9*** For a given fuzzy left derivations subset S , let $\mu_\beta$ be the strongest fuzzy left derivations relation on S , then for $t \in [0,1]$ , we have $(\mu_\beta)_t = \beta_t \times \beta_t$ .

***Proposition 5.10*** For a given fuzzy subset β of KU - algebra X , let $\mu_\beta$ be the strongest left fuzzy derivations relation on X . If $\mu_\beta$ is a fuzzy left derivations KU - ideal of $X \times X$ , then $\beta(d(x)) \leq \beta(d(0)) = \beta(0)$ for all $x \in X$ .

Proof . Since $\mu_\beta$ is a fuzzy left derivations KU- ideal of $X \times X$ , it follows from ($F_1$) that

$\mu_\beta(x, x) = \min \{\beta(d(x)), \beta(d(x))\} \leq \beta(d(0, 0)) = \min \{\beta(d(0)), \beta(d(0))\}$ ,

where $(0, 0) \in X \times X$ then $\beta(d(x)) \leq \beta(d(0)) = \beta(0)$ .

***Remark 5.11*** Let X and Y be KU- algebras , we define $*$ on $X \times Y$ by :

For every $(x, y), (u, v) \in X \times Y$ , $(x, y) * (u, v) = (x * u, y * v)$ , then clearly

$(X \times Y, *, (0, 0))$ is a KU- algebra .

***Theorem 5.12*** Let μ and β be a fuzzy left derivations KU- ideals of KU - algebra X ,then $\mu \times \beta$ is a fuzzy left derivations KU-ideal of $X \times X$ .

Proof : for any $(x, y) \in X \times X$ ,we have ,

$(\mu \times \beta)(d((0, 0))) = \min \{\mu(d(0)), \beta(d(0))\} = \min \{\mu(0), \beta(0)\} \geq$

$\min \{\mu(d(x)), \beta(d(x))\} = (\mu \times \beta)(d(x, y))$ .

Now let $(x_1, x_2), (y_1, y_2), (z_1, z_2) \in X \times X$ , then ,

$(\mu \times \beta)(d(x_1 * z_1, x_2 * z_2)) = \min \{\mu(d(x_1, z_1)), \beta(d(x_2, z_2))\}$

$\geq \min\{\min \{\mu(d(x_1) * (y_1 * z_1))), \mu(d(y_1))\}\ ,\ \min \{\beta(d(x_2) * (y_2 * z_2))), \beta(d(y_2))\}\}$

$$= \min\{\min\{\mu(d(x_1) * (y_1 * z_1))), \mu(d(x_2) * (y_2 * z_2)))\}, \min\{\mu(d(y_1)), \beta(d(y_2))\}\}$$

$$= \min\{(\mu \times \beta)((d(x_1) * (y_1 * z_1), d(x_2) * (y_2 * z_2))), (\mu \times \beta)(d(y_1), d(y_2))\}.$$

Hence $\mu \times \beta$ is a fuzzy left derivations KU- ideal of $X \times X$.

Analogous to theorem 3.2 [16], we have a similar results for fuzzy left derivations KU-ideal, which can be proved in similar manner, we state the results without proof.

***Theorem 5.13*** Let $\mu$ and $\beta$ be a fuzzy left derivations subset of KU-algebra $X$, such that

$\mu \times \beta$ is a fuzzy left derivations KU-ideal of $X \times X$, then

(i) Either $\mu(d(x)) \leq \mu(d(0))$ or $\beta(d(x)) \leq \beta(d(0))$ for all $x \in X$,

(ii) If $\mu(d(x)) \leq \mu(d(0))$ for all $x \in X$, then either $\mu(d(x)) \leq \beta(d(0))$ or $\beta(d(x)) \leq \beta(d(0))$,

(iii) If $\beta(d(x)) \leq \beta(d(0))$ for all $x \in X$, then either $\mu(d(x)) \leq \mu(d(0))$ or $\beta(d(x)) \leq \mu(d(0))$,

(iv) Either $\mu$ or $\beta$ is a fuzzy left derivations KU- ideal of $X$.

***Theorem 5.14*** Let $\beta$ be a fuzzy subset of KU- algebra $X$ and let $\mu_\beta$ be the strongest fuzzy left derivations relation on $X$, then $\beta$ is a fuzzy left derivations KU - ideal of $X$ if and only if $\mu_\beta$ is a fuzzy left derivations KU- ideal of $X \times X$.

proof : Assume that $\beta$ is a fuzzy left derivations KU- ideal $X$, we note from $(F_1)$ that :

$\mu_\beta(0, 0) = \min\{\beta(d(0)), \beta(d(0))\} = \min\{\beta(0), \beta(0)\} \geq \min\{\beta(d(x)), \beta(d(y))\}$

$\qquad = \mu_\beta(d(x), d(y))$.

Now, for any $(x_1,x_2), (y_1,y_2), (z_1,z_2) \in X \times X$, we have from $(F_2)$:

$\mu_\beta(d(x_1 * z_1), d(x_2 * z_2)) = \min\{\beta(d(x_1 * z_1)), \beta(d(x_2 * z_2))\}$

$\qquad \geq \min\{\min\{\beta(d(x_1) * (y_1 * z_1)), \beta(d(y_1))\}, \min\{\beta(d(x_2) * (y_2 * z_2)), \beta(y_2)\}\}$

$\qquad = \min\{\min\{\beta(d(x_1) * (y_1 * z_1)), \beta(d(x_2) * (y_2 * z_2))\}, \min\{\beta(d(y_1)), \beta(d(y_2))\}\}$

$\qquad = \min\{\mu_\beta(d(x_1) * (y_1 * z_1), d(x_2) * (y_2 * z_2)), \mu_\beta(d(y_1), d(y_2))\}.$

Hence $\mu_\beta$ is a fuzzy left derivations KU - ideal of $X \times X$.

Conversely. For all $(x, y) \in X \times X$, we have

Min $\{\beta(0), \beta(0)\} = \mu_\beta(x, y) = \min\{\beta(x), \beta(y)\}$. It follows that

$\beta(0) \geq \beta(x)$ for all $x \in X$, which prove ($F_1$).

Now, let $(x_1, x_2), (y_1, y_2), (z_1, z_2) \in X \times X$, then

$\min\{\beta(d(x_1 * z_1)), \beta d(x_2 * z_2))\} = \mu_\beta(d(x_1 * z_1), d(x_2 * z_2))$

$\geq \min\{\mu_\beta(d(x_1, x_2) * ((y_1, y_2) * (z_1, z_2))), \mu_\beta(d(y_1), d(y_2))\}$

$= \min\{\mu_\beta(d(x_1) * (y_1 * z_1), d(x_2) * (y_2 * z_2)), \mu_\beta(d(y_1), d(y_2))\}$

$= \min\{\min\{\beta(d(x_1) * (y_1 * z_1)), \beta(d(x_2) * (y_2 * z_2))\}, \min\{\beta(d(y_1)), \beta(d(y_2))\}\}$

$= \min\{\min\{\beta(d(x_1) * (y_1 * z_1)), \beta(d(y_1))\}, \min\{\beta((dx_2) * (y_2 * z_2)), \beta(d(y_2))\}\}$

In particular, if we take $x_2 = y_2 = z_2 = 0$, then,

$\beta(d(x_1 * z_1)) \geq \min\{\beta(d(x_1) * (y_1 * z_1)), \beta(d(y_1))\}$ This prove ($FL_2$) and completes the proof.

## *Conclusion*

Derivation is a very interesting and important area of research in the theory of algebraic structures in mathematics. In the present paper, the notion of fuzzy left derivations KU - ideal in KU-algebra are introduced and investigated the useful properties of fuzzy left derivations KU - ideals in KU-algebras.

In our opinion, these definitions and main results can be similarly extended to some other algebraic systems such as BCI-algebra, BCH-algebra ,Hilbert algebra ,BF-algebra -J-algebra ,WS-algebra ,CI-algebra, SU-algebra ,BCL-algebra ,BP-algebra ,Coxeter algebra ,BO-algebra ,PU- algebras and so forth.

The main purpose of our future work is to investigate:

(1) The interval value, bipolar and intuitionistic fuzzy left derivations KU - ideal in KU-algebra.

(2) To consider the cubic structure of left derivations KU - ideal in KU-algebra.

We hope the fuzzy left derivations KU - ideals in KU-algebras, have applications in different branches of theoretical physics and computer science.

## *Algorithm for KU-algebras*

Input ( $X$ : set, $*$ : binary operation)

Output (" $X$ is a KU-algebra or not")

Begin

If $X = \phi$ then go to (1.);

End If

If $0 \notin X$ then go to (1.);

End If

Stop: =false;

$i := 1$;

While $i \leq |X|$ and not (Stop) do

If $x_i * x_i \neq 0$ then

Stop: = true;

End If

$j := 1$

While $j \leq |X|$ and not (Stop) do

If $((y_j * x_i) * x_i) \neq 0$ then

Stop: = true;

End If

End If

$k := 1$

While $k \leq |X|$ and not (Stop) do

If $(x_i * y_j) * ((y_j * z_k) * (x_i * z_k)) \neq 0$ then

Stop: = true;

   End If

  End While

 End While

End While

If Stop then

(1.) Output ("$X$ is not a KU-algebra")

  Else

    Output ("$X$ is a KU-algebra")

  End If

  End.

[26] L.A.Zadeh , Fuzzy sets , inform . and control ,8(1965) , 338-353 .

[27] J. Zhan and Y. L. Liu, "On f-derivations of BCI-algebras," International Journal of Mathematics and Mathematical Sciences, no. 11, pp. 1675–1684, 2005.

Samy M. Mostafa   samymostafa@yahoo.com
Department of Mathematics, Faculty of Education, Ain Shams University, Roxy, Cairo, Egypt.

Ahmed Abd-eldayem   ahmedabdeldayem88@yahoo.com

Department of Mathematics, Faculty of Education, Ain Shams University, Roxy, Cairo, Egypt.